\documentstyle[bezier]{article}

\makeatletter
  
  \@addtoreset{equation}{section}
\makeatother

\newtheorem{theorem}{Theorem}[section]
\newtheorem{proposition}{Proposition}[section]
\newtheorem{conjecture}{Conjecture}[section]
\newtheorem{lemma}{Lemma}[section]

\newtheorem{remark}{Remark}[section]

\newfont{\germ}{eufm10}
\newfont{\germlarge}{eufm10}
\newfont{\slsmall}{cmsl8}

\def\B{{\cal B}}
\def\bbar{\overline{b}}
\def\cd{\cdots}
\def\et#1{\tilde{e}_{#1}}
\def\ft#1{\tilde{f}_{#1}}
\def\geh{\goth{g}}
\def\goth#1{\mbox{\germ #1}}
\def\L{{\cal L}}
\def\La{\Lambda}
\def\la{\lambda}
\def\ot{\otimes}
\def\P{{\cal P}}
\def\pbar{\overline{p}}
\def\Pcl{P_{cl}}
\def\Pcll{(P_{cl}^+)_l}
\def\Proof{\noindent{\sl Proof.}\quad}
\def\Q{{\bf Q}}
\def\qed{~\rule{1mm}{2.5mm}}

\def\Uq{U_q(\geh)}
\def\veps{\varepsilon}
\def\vphi{\varphi}
\def\wt{\mbox{\sl wt}\,}
\def\wts{\mbox{\slsmall wt}\,}
\def\Z{{\bf Z}}

\begin{document}

\title{ Characters of Demazure Modules  \\
        and Solvable Lattice Models  \\ }

\author{Atsuo Kuniba$^1$, Kailash C. Misra$^2$, 
        Masato Okado$^3$,\\
        Taichiro Takagi$^4$ and Jun Uchiyama$^5$} 

\date{ \it $^1$Institute of Physics,\\
       \it University of Tokyo,\\
       \it Komaba, Tokyo 153, Japan\\
       \vskip 4mm
       \it $^2$Department of Mathematics,\\
       \it North Carolina State University,\\
       \it Raleigh, NC 27695-8205, USA\\ 
       \vskip 4mm
       \it $^3$Department of Mathematical Science,\\
       \it Faculty of Engineering Science, \\
       \it Osaka University, Toyonaka, Osaka 560, Japan\\
       \vskip 4mm
       \it $^4$Department of Mathematics and Physics,\\
       \it National Defense Academy, Yokosuka 239, Japan\\
       \vskip 4mm
       \it $^5$Department of Physics,\\
       \it Rikkyo University,\\
       \it Nishi-Ikebukuro, Tokyo 171, Japan}

\maketitle

\begin{abstract}
\noindent
We study the path realization of Demazure crystals
related to solvable lattice models
in statistical mechanics.
Various characters are
represented in a unified way as 
the sums over one dimensional configurations
which we call 
unrestricted, classically restricted and restricted paths.
As an application 
characters of Demazure modules
are obtained 
in terms of $q$-multinomial coefficients
for several level 1 modules of classical affine algebras.
\end{abstract} 

\setcounter{section}{-1}

\section{Introduction}

Let $\Uq$ be a quantum affine algebra and 
$V(\la)$ be the integrable $\Uq$-module with 
highest weight $\la \in P^+$.
Given a Weyl group element $w$, the associated Demazure
module $V_w(\la)$ is the finite dimensional subspace of 
$V(\la)$ generated from the extermal weight space
$V(\la)_{w \la}$ by the $e_i$ generators.
In \cite{Ka}, Kashiwara introduced its crystal
$\B_w(\la)$, which is a finite subset of 
the crystal $\B(\la)$ for $V(\la)$.
With this aid  the 
character of the Demazure module $V_w(\la)$ is expressed as
\begin{equation}
ch V_w(\la) = \sum_{p \in \B_w(\la)} e^{\wts p}.
\label{char}
\end{equation}
The subject of this paper is to calculate (\ref{char}) systematically
by using the path realization of the Demazure crystal
$\B_w(\la)$ studied in \cite{KMOU}, \cite{KMOTU1}.
The realization is based on the earlier one
for $\B(\la)$ \cite{KMN1}, \cite{KMN2} and
has an origin in the analyses of 
solvable lattice models \cite{ABF}, \cite{DJKMO1}, \cite{DJKMO2}.
In these works the object called one dimensional
configuration sums (1dsums) played an
essential role and were studied extensively either 
in their ``infinite lattice limits $j \rightarrow \infty$'' or
``finite truncations $j < \infty$''.
In this paper we consider three kinds of 
1dsums $g_j, \overline{X}_j$ and $X_j$, which we call 
unrestricted, classically restricted and restricted 1dsums, 
respectively.
The unrestricted $g_j$ is relevant to vertex models and so is $X_j$ to 
the restricted solid-on-solid (RSOS) type models.
In section 2 we apply 
the main theorem in \cite{KMOU}
to relate the Demazure character with the 
1dsum $g_j$ for finite $j$.
We shall also clarify the 
relations among the three kinds of 1dsums and thereby 
give a unified picture to understand the 
Demazure characters and various branching functions.
(See Table 1.)
This enables us to evaluate 
these quantities explicitly 
from several known results on the 1dsums.
As an application, in section 3
we shall give $q$-multinomial formulae for $ch V_w(\la)$
for many level 1 modules
$V(\la)$ over $\Uq$ for $\geh$ of classical types
$\geh = A^{(1)}_n, B^{(1)}_n, D^{(1)}_n, 
A^{(2)}_{2n-1}, A^{(2)}_{2n}$ and $D^{(2)}_{n+1}$.
We hope to report on higher level cases in a future
publication.

\section{Path realization of Demazure crystals}

\subsection{Perfect crystals}

Let us recall relevant facts and notations from
\cite{KMN1},\cite{KMN2},\cite{KMOU},\cite{KMOTU1}.
$\alpha_i, \, h_i, \, \Lambda_i\, (i \in I)$ are
the simple roots, coroots and fundamental weights, respectively.
We put $\rho = \sum_{i \in I} \Lambda_i$ and let $\delta$ denote 
the null root.
$P = \oplus_i {\bf Z} \Lambda_i \oplus {\bf Z} \delta$ and
$\Pcl =  \oplus_i {\bf Z} \Lambda_i \subset P$ are 
the weight and the classical weight lattices, respectively.
Set further
$P^+ (\Pcl^+) =\{\la\in P (\Pcl) \mid \langle\la,h_i\rangle\ge0
\mbox{ for any }i\}$ and 
$(\Pcl^+)_l=\{\la\in\Pcl^+\mid \langle\la,c\rangle=l\}$, 
where $c$ is the canonical central element.
For $\la \in P^+$ we let $(\L(\la),\B(\la))$ denote the crystal base
of the irreducible $\Uq$-module $V(\la)$
with highest weight $\la$.
For a crystal base of a finite dimensional 
$\Uq$-module we use the symbol $(L,B)$.
Let $B$ be a perfect crystal of level $l$.
See Definition 4.6.1 in \cite{KMN1} for its definition.
Then for any $\la\in(\Pcl^+)_l$, there uniquely 
exists $b(\la) \in B$ such that 
$\vphi(b(\la)) = \la$.
Here we recall that
$\veps_i(b) = \hbox{max }\{ k \mid \et{i}^k b \neq 0 \}, 
\vphi_i(b) = \hbox{max }\{ k \mid \ft{i}^k b \neq 0 \}, 
\veps(b) = \sum_{i \in I} \veps_i(b) \Lambda_i$ and 
$\vphi(b) = \sum_{i \in I} \vphi_i(b) \Lambda_i$.
Let $\sigma$ be the
automorphism of $(\Pcl^+)_l$ given by $\sigma\la=\veps(b(\la))$.
We put 
$\bbar_k=b(\sigma^{k-1}\la)$ and $\la_k=\sigma^k\la$. 
Then perfectness
assures that we have the isomorphism of crystals
\begin{equation}
\B(\la_{k-1})\simeq \B(\la_k)\ot B.
\label{iso1}
\end{equation}
Define the set of paths $\P(\la,B)$ by
\[
\P(\la,B)=\{p=\cdots\ot p(2)\ot p(1)\mid p(j)\in B,p(k)=
\bbar_k\mbox{ for }k\gg0\},
\]
By iterating (\ref{iso1}) we have an isomorphism of crystals
\begin{equation}
\B(\la) \simeq \P(\la,B).
\label{iso2}
\end{equation}
In particular, the image of 
the highest weight vector $u_\la \in  \B(\la)$ is given by 
$\pbar = \cd\ot\bbar_k\ot\cd\ot\bbar_2\ot
\bbar_1$.
We call $\pbar$ the {\em ground-state} path.
The actions of 
$\et{i}$ and $\ft{i}$ on $\P(\la,B)$ are determined 
explicitly by the signature rule.
See section 1.3 of \cite{KMOU}.

To describe the weights 
on $\P(\la, B)$ it is necessary to introduce the 
energy function $ H:B\ot B\rightarrow\Z$.
Up to an additive constant it is determined 
by requiring the following 
for any $b,b'\in B$ and $i\in I$ such that $\et{i}(b\ot b')\ne0$.
\begin{equation}
H(\et{i}(b\ot b'))=\left\{\begin{array}{ll}
H(b\ot b')&\quad\mbox{if }i\ne0\\
H(b\ot b')+1&\quad\mbox{if }i=0\mbox{ and }\vphi_0(b)\ge\veps_0(b')\\
H(b\ot b')-1&\quad\mbox{if }i=0\mbox{ and }\vphi_0(b)<\veps_0(b').
\end{array}\right.
\label{defh}
\end{equation}
Under the 
isomorphism (\ref{iso2}), the weight of a path
$p=\cdots\ot p(2)\ot p(1)$ is given by 
(Proposition 4.5.4 in \cite{KMN1})
\begin{eqnarray}
\wt p&=&\la+\sum_{i=1}^\infty
\left(\wt p(i) - \wt \bbar_i \right)\nonumber\\
&&\quad-\left(\sum_{i=1}^\infty 
i(H(p(i+1)\ot p(i))-H(\bbar_{i+1}\ot \bbar_i))\right)\delta. 
\label{eq:weight}
\end{eqnarray}
We remark the weight relation
\begin{equation}
\la_j = \la - \sum_{i=1}^j \wt \bbar_i,
\label{weightrelation}
\end{equation}
which is valid for any $j \ge 0$.

\subsection{Demazure modules}

Let $\{r_i\}_{i\in I}$ be the set of simple reflections, and let
$W$ be the Weyl group. 
For $\la \in (\Pcl^+)_l$ we consider the
Demazure module $V_w(\la)$ generated from the 
extremal weight space $V(\la)_{w\la}$. 
By definition its character is given by
$ch V_w(\la) = \sum_\mu {\rm dim }(V_w(\la))_\mu e^\mu$.
For $\mu \in P, \, i \in I$ define the operator
$D_i: {\bf Z}[P] \rightarrow {\bf Z}[P]$ by
$$
D_i(e^\mu) = {e^{\mu + \rho} - e^{r_i(\mu + \rho)} \over
1 - e^{-\alpha_i}}e^{-\rho}.
$$
Let $w = r_{i_k} \cdots r_{i_1} \in W$ be a reduced expression.
Then the following character formula is well known 
\cite{Dem}, \cite{Ku}, \cite{Mat}.
\begin{equation}
ch V_w(\la) = D_{i_k} \cdots D_{i_2} D_{i_1}(e^\mu).
\label{kumar}
\end{equation}
{}From this one has a recursion relation
\begin{equation}
ch V_{r_i w}(\la) = D_i\left( ch V_w(\la) \right) 
\quad{\rm if }\, r_i w \succ w.
\label{eq:demazurerec}
\end{equation}

Let $(\L(\la),\B(\la))$ be the crystal base of $V(\la)$. 
In \cite{Ka}
Kashiwara showed that for each $w\in W$, there exists a subset
$\B_w(\la)$ of $\B(\la)$ such that
\begin{equation}
\frac{V_w(\la)\cap\L(\la)}{V_w(\la)\cap q\L(\la)}
=\bigoplus_{b\in\B_w(\la)}\Q b.
\end{equation}
Furthermore, $\B_w(\la)$ has the following recursive property.
\begin{eqnarray}
&&\mbox{If }r_iw\succ w,\mbox{ then}\nonumber\\
&&\B_{r_iw}(\la)=\bigcup_{n\ge0}\ft{i}^n\B_w(\la)\setminus\{0\},
\label{recursive}
\end{eqnarray}
which is analogous to (\ref{eq:demazurerec}).
We call $\B_w(\la)$ a {\em Demazure crystal}.
One can now express $ch V_w(\la)$ 
as in  (\ref{char}).
It affords an effecient way to calculate
the Demazure character through
the path realization of $\B_w(\la)$ given in the next subsection
and (\ref{eq:weight}).
Relations with the formula 
(\ref{kumar}) and the 1dsums 
will also be explained in section 2.

\subsection{Path realization}

In \cite{KMOU} the image of $\B_w(\la)$ under the 
isomorphism (\ref{iso2}) is determined for a 
suitably chosen Weyl group element $w$.
Let us recall the main theorem therein, which gives a 
path realization of the Demazure crystal.
There appears the {\em mixing index} $\kappa$ 
specified from $\la$ and $B$.
(See section 2.3 of \cite{KMOU}.)
In this paper we shall only consider the case $\kappa=1$.
Let $\la$ be an element of $\Pcll$, and let $B$ be a classical crystal. 
For the theorem, we need to assume four conditions (I-IV).
\begin{description}
\item[(I)   ]
$B$ is perfect of level $l$.
\end{description}
Thus, we can assume an isomorphism between $\B(\la)$ and the set of paths
$\P(\la,B)$. Let $\pbar=\cd\ot\bbar_2\ot\bbar_1$ denote the ground state 
path. Fix a positive integer $d$. For a set of elements
$i_a^{(j)}$ ($j\ge1,1\le a\le d$) in $I$, we define
$B_a^{(j)}$ ($j\ge1,0\le a\le d$) by 
\[
B^{(j)}_0=\{\bbar_j\},\hspace{1cm}
B_a^{(j)}=\bigcup_{n\ge0}\ft{i_a^{(j)}}^n B_{a-1}^{(j)}\setminus\{0\}
\quad(a=1,\cdots,d).
\]
\begin{description}
\item[(II)  ]
For any $j\ge1$,
$B_d^{(j)}=B$.
\item[(III) ]
For any $j\ge1$ and $1\le a\le d$,
$\langle\la_j,h_{i^{(j)}_a}\rangle\le\veps_{i^{(j)}_a}(b)$
for all $b\in B^{(j)}_{a-1}$.
\end{description}
We now define an element $w^{(k)}$ of the Weyl group $W$ by
\[
w^{(0)}=1,\hspace{1cm}
w^{(k)}=r_{i^{(j)}_a}w^{(k-1)}\quad\mbox{for }k>0,
\]
where $j$ and $a$ are fixed from $k$ by $k=(j-1)d+a,j\ge1,1\le a\le d$.
\begin{description}
\item[(IV)  ]
$w^{(0)}\prec w^{(1)}\prec\cd\prec w^{(k)}\prec\cd$.
\end{description}
See \cite{KMOU}, \cite{KMOTU1} on how to check the last condition.

Finally we define a subset $\P^{(k)}(\la,B)$ of $\P(\la,B)$ as follows.
We set $\P^{(0)}(\la,B)=\{\overline{p}\}$.
For $k>0$, 
\begin{equation} \label{eq:def_P}
\P^{(k)}(\la,B) = 
\cdots\ot B^{(j+2)}_0\ot B^{(j+1)}_0\ot B^{(j)}_a
\ot B^{\ot(j-1)},
\end{equation}
where $j \ge 1$ and $1 \le a \le d$ are 
uniquely specified by $k=(j-1)d+a$.

Now we have 

\begin{theorem}[{\rm \cite{KMOU}}] \label{iso3}
Under the assumptions (I-IV), we have 
\[
\B_{w^{(k)}}(\la)\simeq\P^{(k)}(\la,B).
\]
\end{theorem}

The proof is done by showing the recursion relation 
(\ref{recursive}) in the path realization.

\section{One dimensional sums}

Here we first introduce the unrestricted 1 dimensional
sum (1dsum) $g_j$ 
and express the Demazure characters in terms of it.
After establishing its fundamental properties 
we then introduce 
classically restricted 1dsums $\overline{X}_j$ and restricted 
1dsums $X_j$ and study their relations.
In the working below we shall use the variable
\[
q = e^{-\delta}.
\]

\subsection{Unrestricted 1dsum}

For $j \in {\bf Z}_{\ge 0}, \, b \in B$ and 
$\mu \in P$, put 
\begin{equation}
\P_j(b, \mu) = \{ b \ot b_j \ot \cdots \ot b_1 \in B^{\ot (j+1)} \vert
wt(b_j \ot \cdots \ot b_1) = cl(\mu) \}.
\label{defpj}
\end{equation}
In the sequel we always assume the relation 
$k = (j-1)d + a$.
Comparing (\ref{defpj}) with (\ref{eq:def_P}) we have
$$
\P^{(k)}(\la, B) = \sqcup_{\mu \in P_{cl}}
\sqcup_{b \in B^{(j)}_a} 
\P_{j-1}(b,\mu).
$$
For $j \in {\bf Z}_{\ge 0}, \, b \in B$ and 
$\mu \in P$ we define the (unrestricted) 1dsum
as follows.
\begin{equation}
g_j(b,\mu) =  q^{(\Lambda_0, \mu)} 
\sum_{ b_{j+1} \ot \cdots \otimes b_1 \in \P_j(b,\mu)} 
q^{\sum_{i=1}^jiH(b_{i+1} \ot b_i)}.
\label{def1dsum}
\end{equation}
Defining a map $E: B^{\ot (j+1)} \rightarrow {\bf Z}$ by
\begin{equation}
E(b_{j+1}\ot \cdots \ot b_1) = \sum_{i=1}^j i H(b_{i+1}\ot b_i),
\label{eq:energy}
\end{equation}
one can write (\ref{def1dsum}) as
\begin{equation}
g_j(b, \mu) = q^{(\Lambda_0, \mu)}
\sum_{p \in \P_j(b,\mu)} q^{E(p)}.
\label{eq:def2}
\end{equation}
Note that
\begin{eqnarray}
&&g_j(b, \mu) = 0 \, \, \mbox{ unless } \langle\mu ,c\rangle = 0 \nonumber \\
&&g_j(b,\mu+m\delta) = q^m g_j(b,\mu) \,\, \mbox{ for any } \mu \in P
\mbox{ and } m \in {\bf Z}. \nonumber \\
\label{1dsumproperty}
\end{eqnarray}
For the Weyl group element $w^{(k)}$ in ({\bf III}),
the Demazure character (\ref{char}) is expressed in terms of
the 1dsum.
\begin{proposition} \label{pr:by1dsum}
\begin{eqnarray}
&&ch V_{w^{(k)}}(\la) = q^{-c_j}
\sum_{\mu \in P_{cl}} e^{\la_j + \mu} 
\sum_{b \in B^{(j)}_a} q^{jH(\bbar_{j+1} \ot b)} g_{j-1}(b, \mu-wt(b)),
\nonumber \\
&&c_j = \sum_{i=1}^j i H(\bbar_{i+1} \ot \bbar_i).
\label{by1dsum}
\end{eqnarray}
\end{proposition}

\Proof
Substitute (\ref{weightrelation}) and (\ref{def1dsum})
into the rhs.
Setting 
$p = \cdots \ot \bbar_{j+1} \ot b_j \ot \cdots \ot b_1$
and noting (\ref{eq:weight}) one finds that the result
is equal to 
$\sum_{p \in \P^{(k)}(\la, B)} e^{\wts p}$.
Thus the assertion follows from Theorem \ref{iso3}.
\qed

The 1dsums are uniquely characterized also from the 
recursion relation and the initial condition as follows.

\begin{proposition} \label{pr:rec1dsum}
\begin{eqnarray}
&& g_j(b,\mu) = \sum_{b^\prime \in B} q^{j H(b \ot b^\prime)} 
g_{j-1}(b^\prime,\mu-wt(b^\prime)) \nonumber\\
&& g_0(b,\mu) = \delta_{0 \mu}.
\label{rec1dsum}
\end{eqnarray}
\end{proposition}

\Proof
In the definition (\ref{def1dsum}) 
put $b^\prime = b_j$ and notice 
$b_{j+1} = b$ due to (\ref{defpj}).
\qed

When $k =jd \, (a=d), B^{(j)}_a = B$ due to
{\bf (II)}.
Thus Proposition \ref{pr:rec1dsum} simplifies the sum
in (\ref{by1dsum}) into 
\begin{equation}\label{simplecase}
ch V_{w^{(jd)}}(\la) = q^{-c_j}
\sum_{\mu \in P_{cl}} e^{\la_j + \mu} 
g_j(\bbar_{j+1}, \mu).
\end{equation}

The following relation is of primary importance in later
discussion.
See Remark \ref{fundamental} and the proof of Proposition
\ref{pr:xbyg}.

\begin{proposition} \label{pr:2mrelation}
For $b \in B$, let $m = \vphi_i(b)$. Then we have
\begin{equation}
\sum_{t=0}^m g_j(\ft{i}^t b,\mu+t\alpha_i)q^{tj\delta_{0 i}}=
\sum_{t=0}^m g_j\bigl(\ft{i}^t b,r_i(\mu+(m-t)\alpha_i)\bigr)
q^{tj\delta_{0 i}}.
\label{eq:2mrelation}
\end{equation}
\end{proposition}

For the proof we need a few Lemmas.
The following is an immediate consequence of the signature rule.

\begin{lemma} \label{le:inequality}
For any $b_1, b_2 \in B$ and $i \in I$ we have
\begin{eqnarray}
&&\vphi_i(b_1\ot b_2) \ge \vphi_i(b_1) + \langle h_i, \wt b_2 \rangle
\nonumber \\
&&\veps_i(b_1\ot b_2) \ge - \langle h_i, \wt(b_1\ot  b_2) \rangle.
\end{eqnarray}
\end{lemma}

\begin{lemma} \label{le:bijection}
Let $n = \langle h_i, \mu \rangle + m$ and assume $n \ge 0$.
Then the map
\begin{equation}
\ft{i}^n : \sqcup_{t=0}^m 
\P_j(\ft{i}^t b, \mu + t\alpha_i) \rightarrow
\sqcup_{t=0}^m 
\P_j\left(\ft{i}^t b, r_i(\mu + (m-t)\alpha_i)\right)
\label{eq:bijection}
\end{equation}
is a bijection.
\end{lemma}

\Proof
The image is certainly within the rhs by the weight reason
unless it is zero.
Since $\ft{i}^n p = p^\prime$ is equivalent to 
$\et{i}^n p^\prime = p$, it suffices to show for 
$0 \le t \le m$ that 
$\vphi_i(p) \ge n$ for any 
$p \in \P_j(\ft{i}^t b, \mu + t\alpha_i)$ and
$\veps_i(p^\prime) \ge n$ for any
$p^\prime \in \P_j\left(\ft{i}^t b, r_i(\mu + (m-t)\alpha_i)\right)$.
Applying Lemma \ref{le:inequality}, one has
\begin{eqnarray*}
&\vphi_i(p) & \ge m-t + \langle h_i, cl(\mu + t\alpha_i) \rangle = n + t, \\
&\veps_i(p^\prime) & \ge 
\veps_i(\ft{i}^t b) - \vphi_i(\ft{i}^t b) - 
\langle h_i, cl\left(r_i(\mu + (m-t)\alpha_i)\right) \rangle \\
&& \ge 2t-m + \langle h_i, \mu + (m-t)\alpha_i \rangle = n,
\end{eqnarray*}
which completes the proof.
\qed

\begin{lemma} \label{le:energy}
Let $b \in B, \, m = \vphi_i(b), \, \xi \in P$.
For $0 \le t \le s \le m$, assume that 
$p \in \P_j(\ft{i}^s b, \xi)$ and
$\et{i}^n p \in \P_j(\ft{i}^t b, \xi+(n+t-s)\alpha_i)$.
Then we have
$$
E(\et{i}^n p) = E(p) + \left( (j+1)(s-t) - n \right)\delta_{0 i}.
$$
\end{lemma}

\Proof
Use (\ref{defh}) in (\ref{eq:energy}).
\qed
\vskip0.2cm
\noindent
{\em Proof of Proposition \ref{pr:2mrelation}}.
Put $n = \langle h_i, \mu \rangle + m$.
Since (\ref{eq:2mrelation}) is invariant under
the change $\mu \rightarrow r_i(\mu + m\alpha_i)$,
we may assume $n \ge 0$ with no loss of generality.
The lhs of (\ref{eq:2mrelation}) is written as
\begin{equation}
\sum_{0 \le t \le s \le m} 
q^{(\Lambda_0, \mu + t\alpha_i)+ tj\delta_{0 i}}
\sum^{}_{p^\prime}{}^\prime q^{E(p^\prime)},
\label{eq:1}
\end{equation}
where $\sum^{}_{p^\prime}{}^\prime$ extends over those
$p^\prime \in \P_j(\ft{i}^t b, \mu + t\alpha_i)$ such that
$\ft{i}^n p^\prime \in \P_j(\ft{i}^s b, r_i(\mu + (m-s)\alpha_i))$.
The sum 
$\sum_{0 \le t \le s \le m} 
\sum^{}_{p^\prime}{}^\prime$ in total 
ranges over the set appearing in the lhs of 
(\ref{eq:bijection}).
Thus Lemma \ref{le:bijection} allows a change of the summation
variable into $p = \ft{i}^n p^\prime$, and thereby transforms
(\ref{eq:1}) into
\begin{equation}
\sum_{0 \le t \le s \le m} 
q^{(\Lambda_0, \mu + t\alpha_i)+ tj\delta_{0 i}}
\sum^{}_p{}^{\prime\prime} q^{E(\et{i}^n p)}.
\label{eq:2}
\end{equation}
Here 
$\sum^{}_{p}{}^{\prime\prime}$ extends over those
$p \in \P_j(\ft{i}^s b, r_i(\mu + (m-s)\alpha_i))$
such that 
$\et{i}^n p \in \P_j(\ft{i}^t b, \mu + t\alpha_i)$.
Setting $\xi = r_i(\mu + (m-s)\alpha_i)$,
one can apply Lemma \ref{le:energy} to rewrite (\ref{eq:2}) as
\begin{eqnarray}
&&\sum_{0 \le t \le s \le m} 
q^{(\Lambda_0, r_i(\mu + (m-s)\alpha_i))+ sj\delta_{0 i}}
\sum^{}_p{}^{\prime\prime} q^{E(p)}\nonumber \\
&&=\sum_{0 \le s \le m} 
q^{(\Lambda_0, r_i(\mu + (m-s)\alpha_i))+ sj\delta_{0 i}}
\sum^{}_{p \in \P_j(\ft{i}^sb, r_i(\mu+(m-s)\alpha_i))}
q^{E(p)}.\nonumber
\end{eqnarray}
By (\ref{eq:def2}) 
the last expression is the rhs of (\ref{eq:2mrelation}).
\qed

The relation (\ref{eq:demazurerec}) for the 
Demazure character (\ref{by1dsum}) implies yet another
recursion relation for the 1dsums than Proposition \ref{pr:rec1dsum}.
\begin{proposition}\label{pr:1dsumdemarec}
For any $0 \le a \le d-1$ and $\mu \in P_{cl}$ one has
\begin{eqnarray*}
&&\sum_{b \in B^{(j)}_{a+1} \setminus B^{(j)}_a}
q^{jH(\bbar_{j+1} \ot b)} g_{j-1}(b, \mu-wt(b)) \nonumber \\
&& =
\sum_{b \in B^{(j)}_{a+1}}
q^{jH(\bbar_{j+1} \ot b)} g_{j-1}(b, \mu+\alpha_{i^{(j)}_{a+1}}-wt(b))
\nonumber \\
&& - \sum_{b \in B^{(j)}_a}
q^{jH(\bbar_{j+1} \ot b)} 
g_{j-1}(b, r_{i^{(j)}_{a+1}}(\mu+\rho+\la_j)-\la_j-\rho-wt(b)).
\end{eqnarray*}
\end{proposition}

\Proof
Suppose $1 \le a \le d-1$.
Substitute (\ref{by1dsum}) into  
(\ref{eq:demazurerec}) with $w=w^{(k)}$ and 
$i = i^{(j)}_{a+1}$.
Note that 
$D_i(e^{\mu + z \delta}) = e^{z \delta} D_i(e^\mu)$
holds for any $z \in {\bf Z}$ and $\mu \in P_{cl}$
because of $r_i(\delta) = \delta$.
After multiplying both sides by $1 - e^{-\alpha_i}$,
comparison of the coefficients of 
$e^{\lambda_j + \mu}$ leads to the above relation.
The case $a=0$ is similar.
\qed
\begin{remark}\label{fundamental}
It is possible to prove Proposition \ref{pr:1dsumdemarec}
without using (\ref{eq:demazurerec}) and 
(\ref{by1dsum}) but only from (\ref{eq:2mrelation}).
In this sense, (\ref{eq:2mrelation}) is the most fundamental
relation of the 1dsums implied from the Demazure 
recursion relation (\ref{eq:demazurerec}).
\end{remark}

\subsection{Classically restricted and restricted 1dsum}

The 1dsums discussed so far is related to vertex models.
Now we proceed to the two variants of them related to 
restricted solid-on-solid (RSOS) type models 
(cf. \cite{ABF}, \cite{DJKMO1}) and Kostka-type polynomials 
(cf.\cite{Ki}, \cite{NY}).
Here they shall be called the restricted 1dsums and 
the classically restricted 1dsums, respectively.
Let $\overline{\Lambda}_i = 
\Lambda_i - \langle c , \Lambda_i \rangle \Lambda_0$ and 
put 
$\overline{P}_{cl}^+ = \oplus_{i \in I \setminus \{0 \}}
{\bf Z}_{\ge 0} \overline{\Lambda}_i$.
Fix a non-negative integer $l'$.
Given $\xi \in (P^+_{cl})_{l + l'}$
(resp. $\overline{\xi} \in \overline{P}_{cl}^+$) and $b \in B$ 
we define
\begin{eqnarray*}
&&(\xi, b) \,\hbox{ is {\em  admissible} } \Leftrightarrow 
\et{i}^{\langle h_i, \xi \rangle + 1} b = 0 \,  \,
\forall i \in I,\\
&&(\overline{\xi}, b) \,\hbox{ is {\em classically  admissible} }
 \Leftrightarrow 
\et{i}^{\langle h_i, \overline{\xi} \rangle + 1} b = 0 \,\, 
\forall i \in I \setminus \{ 0 \}.\\
\end{eqnarray*}
Recall that $l$ is the level of the perfect crystal $B$.
It is easy to see that if $(\xi, b)$ is admissible
(resp. $(\overline{\xi}, b)$ is classically admissible)  then
$\xi + wt(b) \in (P^+_{cl})_{l + l'}$ 
(resp. $\overline{\xi} + wt(b) \in \overline{P}_{cl}^+$).
The admissibility condition has been introduced by \cite{DJO}
in the study of $q$-vertex operators by the crystal base theory.
For $j \in {\bf Z}_{\ge 0}, b \in B$ and 
$\xi, \eta \in (P^+_{cl})_{l + l'}$ such that
$(\xi-wt(b), b)$ is admissible, 
(resp. $\overline{\xi}, \overline{\eta} \in \overline{P}_{cl}^+$ such that 
$(\overline{\xi}-wt(b), b)$ is classically admissible), 
we define $q$-polynomials 
$X_j(b,\xi, \eta)$ and 
$\overline{X}_j(b,\overline{\xi},\overline{\eta})$ to be the sum
$$
\sum_{ b_j, \ldots, b_1 \in B, b_{j+1} = b}
q^{\sum_{i=1}^jiH(b_{i+1} \ot b_i)}.
$$
Here the outer sum $\sum$ is taken over $b_j, \ldots, b_1 \in B$
under the following conditions for each case.
\begin{eqnarray*}
X_j(b, \xi, \eta)\, {\rm case}:
&&\xi_i + wt(b_i) = \xi_{i-1} \, \mbox{ for } 1 \le i \le j, \quad
\xi_j = \xi, \, \xi_0 = \eta,\\
&&(\xi_i, b_i) \mbox{ is admissible for } 1 \le i \le j.\\
\overline{X}_j(b, \overline{\xi}, \overline{\eta})\, {\rm case}:
&&\overline{\xi}_i + wt(b_i) = \overline{\xi}_{i-1} \, 
\mbox{ for } 1 \le i \le j, \quad
\overline{\xi}_j = \overline{\xi}, \, 
\overline{\xi}_0 = \overline{\eta},\\
&&(\overline{\xi_i}, b_i) \mbox{ is classically admissible for } 
1 \le i \le j.\\
\end{eqnarray*}
We define $X_j(b, \xi, \eta)$ 
(resp. $\overline{X}_j(b,\overline{\xi},\overline{\eta})$) to be zero if
$(\xi-wt(b),b)$ is not admissible
(resp. if $(\overline{\xi}-wt(b),b)$ is not classically admissible).
This implies that 
$\varphi_i(b) \le \langle h_i, \xi \rangle$
for all $i \in I$
(resp. $\varphi_i(b) \le \langle h_i, \overline{\xi} \rangle$
for all $i \in I\setminus \{ 0 \}$).
We shall call $X_j(b, \xi, \eta)$
and $\overline{X}_j(b, \overline{\xi}, \overline{\eta})$ 
the (level $l + l'$) restricted 1dsums and 
the classically restricted 1dsums, respectively.

For any $\overline{\xi} \in \overline{P}_{cl}^+$,
$\overline{\xi} + (l + l')\Lambda_0$ belongs to 
$(P^+_{cl})_{l+l'}$ for sufficiently large $l'$.
Moreover the pair $(\overline{\xi},b)$ is classically admissible
if and only if 
$(\overline{\xi} + (l + l')\Lambda_0, b)$ is admissible in the limit 
$l' \rightarrow \infty$.
Therefore we have
\begin{proposition} \label{pr:xandxbar}
For any $\overline{\xi}, \overline{\eta} \in \overline{P}_{cl}^+$,
\begin{eqnarray*}
\overline{X}_j(b, \overline{\xi}, \overline{\eta}) = 
\lim_{l' \rightarrow \infty} X_j(b, \overline{\xi} + (l+l')\Lambda_0,
\overline{\eta} + (l+l')\Lambda_0). \nonumber \\
\end{eqnarray*}
\end{proposition}
As Proposition \ref{pr:rec1dsum}, 
the 1dsums $X_j(b, \xi, \eta)$
and $\overline{X}_j(b, \overline{\xi}, \overline{\eta})$ 
are characterized by the recursion relation and the 
initial condition as follows.
\begin{proposition} \label{pr:recr1dsum}
\begin{eqnarray*}
&& X_j(b, \xi, \eta) = 
\sum_{b^\prime \in B, \, (\xi, b'): \hbox{\scriptsize admissible}} 
q^{j H(b \ot b^\prime)} 
X_{j-1}(b^\prime,\xi + wt(b^\prime), \eta),\nonumber \\
&& X_0(b, \xi, \eta) = \delta_{\xi  \eta}, \nonumber \\
&& \overline{X}_j(b, \overline{\xi}, \overline{\eta}) = 
\sum_{b^\prime \in B, \, (\overline{\xi}, b'): 
\hbox{\scriptsize classically admissible}} 
q^{j H(b \ot b^\prime)} 
\overline{X}_{j-1}
(b^\prime,\overline{\xi} + wt(b^\prime), \overline{\eta}),\nonumber \\
&& \overline{X}_0(b, \overline{\xi}, \overline{\eta}) = 
\delta_{\overline{\xi} \overline{\eta}}.
\end{eqnarray*}
\end{proposition}
In order to express $X_j$ and $\overline{X}_j$ in terms of 
$g_j$, we need to assume
\begin{conjecture}\label{con:disjoint}
For any $\xi \in (P_{cl}^+)_l$, there exists 
a disjoint union decomposition
(not necessarily unique) as
\begin{eqnarray}
&&\{ b \in B \mid \, (\xi,b) \hbox{ is not admissible} \} \nonumber \\
&&= \bigsqcup_{b' \in B, 
\varepsilon_i(b') = \langle h_i, \xi \rangle + 1\, 
\hbox{\scriptsize for some } i \in I}
\{ \ft{i}^t b' \mid 0 \le t \le \varphi_i(b') \}.
\label{affinedisjoint}
\end{eqnarray}
\end{conjecture}
We have proved this for $U_q(A^{(1)}_n), \, B = l$-fold symmetric 
tesnor case and have checked several other cases.
This property seems reflecting an intrinsic combinatorial nature
of perfect crystals.
If the conjecture holds, it follows by the same argument 
as for Proposition \ref{pr:xandxbar} that 
for any $\overline{\xi} \in \overline{P}_{cl}^+$, there exists 
a disjoint union decomposition
(not necessarily unique) as
\begin{eqnarray}
&&\{ b \in B \mid \, (\overline{\xi},b) 
\hbox{ is not classically admissible} \} \nonumber \\
&&= \bigsqcup_{b' \in B, 
\varepsilon_i(b') = \langle h_i, \overline{\xi} \rangle + 1\,
\hbox{\scriptsize for some } i \in I\setminus \{0 \}}
\{ \ft{i}^t b' \mid 0 \le t \le \varphi_i(b') \}.
\label{classicaldisjoint}
\end{eqnarray}
\begin{proposition} \label{pr:xbyg}
If Conjecture \ref{con:disjoint} holds, the
restricted and classically restricted 1dsums are 
expressed as linear superpositions
of the 1dsum $g_j$ over the affine Weyl group $W$ and 
the classical Weyl group $\overline{W}$, respectively  as
\begin{eqnarray}
&&X_j(b, \xi, \eta) = 
\sum_{w \in W} \mbox{ det } w\, 
g_j(b, w(\eta + \rho) - \xi - \rho),
\label{xbyg}\\
&&\overline{X}_j(b, \overline{\xi}, \overline{\eta}) = 
\sum_{w \in \overline{W} } \mbox{ det } w\, 
g_j(b, w(\overline{\eta} + \overline{\rho}) 
- \overline{\xi} - \overline{\rho}).
\label{xbarbyg}
\end{eqnarray}
\end{proposition}

\Proof
We show (\ref{xbyg}) from (\ref{affinedisjoint}).
(\ref{xbarbyg}) can be verified from (\ref{classicaldisjoint})
analogously.
Let $F_j(b, \xi, \eta)$ denote the rhs of (\ref{xbyg}).
We are to show that $F_j(b, \xi, \eta)$ fulfills the 
properties in Proposition \ref{pr:recr1dsum}.
To check the initial condition is easy.
By using (\ref{rec1dsum}) one has
$$
F_j(b, \xi, \eta) = \sum_{b' \in B} q^{jH(b \ot b')}
F_{j-1}(b', \xi + wt(b'), \eta).
$$
The sum here is similar to the one in
Proposition \ref{pr:recr1dsum}
but without the constraint 
$(\xi, b'):\hbox{ admissible}$. 
Thus it is enough to show the cancellation of those
unwanted contributions from non-admissible $b'$, namely,
$$
0 = \sum_{b' \in B, (\xi,b'): \hbox{non-admissible}} 
q^{jH(b \ot b')}
F_{j-1}(b', \xi + wt(b'), \eta).
$$
Under the assumption (\ref{affinedisjoint})
it suffices to show the further decomposed form of this as 
$$
0 = \sum_{t=0}^m 
q^{jH(b \ot \ft{i}^t b')}
F_{j-1}(\ft{i}^t b', \xi + wt(\ft{i}^t b'), \eta)
$$
for each $b' \in B$ such that 
$\varepsilon_i(b') = \langle h_i, \xi \rangle + 1$.
Here $m = \varphi_i(b') = \langle h_i, \xi + \rho + wt(b') \rangle$.
{}From $\varphi_i(b) \le \langle h_i, \xi \rangle$  
one has
$H(b \ot \ft{i}^t b') = H(b \ot b') + t\delta_{0 i}$.
By using 
$wt(\ft{i}^t b') = wt(b') - t \alpha_i + t \delta_{0 i} \delta$ 
and (\ref{1dsumproperty}) further, 
the rhs of the above is expressed as
$$
q^{jH(b \ot b')} \sum_{w \in W} \mbox{ det }w 
\sum_{t=0}^m g_{j-1}(\ft{i}^tb', 
w(\eta+\rho) - \xi - \rho - wt(b') + t\alpha_i) 
q^{t(j-1)\delta_{0 i}}.
$$
Upon applying (\ref{eq:2mrelation}), one finds that 
this quantity is precisely equal to itself with $w(\eta+\rho)$ replaced by 
$r_iw(\eta+\rho)$.
Thus it vanishes because of 
$\mbox{ det } r_i w = - \mbox{ det } w$.
\qed

\subsection{Relation with affine Lie algebra and coset characters}

Given a $\Uq$ module $M$ and $\mu \in P$,
let $[M : \mu ]$ (resp. $[M : \mu ]_{cl}$) denote the dimension
of the linear space 
$\{ v \in M \mid wt(v) = \mu, \,
e_i v = 0 
\hbox{ for all } i \in I (\hbox{resp. } i \in I \setminus \{ 0 \})  \}$.
Let $c_j$ be as in Proposition \ref{pr:by1dsum}, 
$\lambda, \lambda_j \in (P_{cl}^+)_l, 
V(\lambda), V_w(\lambda), \bbar_j, b(\lambda), d$ and $\sigma$ be
as in section 1.
Put $\overline{\lambda}_j = \lambda_j - l \Lambda_0$.
The $j \rightarrow \infty$ limits of 
$g_j, X_j$ and $\overline{X}_j$ give rise to various branching functions.
We summarize them in
\begin{proposition}\label{pr:branching}
For $\mu \in P,
\xi \in (P_{cl}^+)_{l'}, \eta \in (P_{cl}^+)_{l+l'}$ and 
$\overline{\eta} \in \overline{P}_{cl}^+$ we have
\begin{eqnarray}
&&\lim_{j \rightarrow \infty} q^{-c_j} g_j(\bbar_{j+1}, \mu) = 
\sum_i \left( \dim V(\lambda)_{\mu - i \delta} \right) q^i,
\label{glimit}\\
&&\lim_{j \rightarrow \infty} q^{-c_j} 
X_j(\bbar_{j+1}, \xi + \lambda_j, \eta) = 
\sum_i [V(\xi) \otimes V(\lambda) : \eta - i \delta] \, q^i,
\label{xlimit}\\
&&\lim_{j \rightarrow \infty} q^{-c_{j}} 
\overline{X}_{j}(\bbar_{j+1}, \overline{\lambda}_j, \overline{\eta}) = 
\sum_i [V(\lambda) : 
\overline{\eta} + l\Lambda_0 - i \delta]_{cl}\, q^i.
\label{xbarlimit}\\
\nonumber
\end{eqnarray}
\end{proposition}

\Proof
(\ref{glimit}) is due to \cite{KMN1}.
(\ref{xlimit}) is due to
\cite{DJO}. 
To show (\ref{xbarlimit}) recall that 
any path can be written in the form 
$p = u_{\lambda_j} \otimes b_j \otimes \cdots \otimes b_1$
$(b_1, \ldots, b_j \in B)$ for sufficiently large $j$,
where one may identify   
$u_{\lambda_j} = \cdots \otimes \bbar_{j+2} \otimes \bbar_{j+1}$.
For a path $p = u_{\lambda_j} \otimes b_j \otimes \cdots \otimes b_1$,
the condition
$\et{i}p = 0\, \forall i \in I \setminus \{ 0 \}$ is
equivalent to the requirement that
$( \overline{\lambda}_j + wt(b_j) + \cdots + wt(b_{i+1}), b_i)$ is 
classically admissible for $1 \le i \le j$.
Since the weight of the path $p$ is given by 
$\overline{\eta} + l \Lambda_0 - \left( 
E(\bbar_{j+1} \otimes b_j \otimes \cdots \otimes b_1) - c_j \right) \delta$,
this completes the proof of (\ref{xbarlimit}).
\qed

Up to an overall power of $q$, (\ref{glimit}) is a string function 
\cite{KP}, (\ref{xlimit}) is a branching coefficient of the module
$V(\eta)$ in the tensor product 
$V(\xi) \otimes V(\lambda)$, 
(\ref{xbarlimit}) is the 
branching coefficient of 
the irreducible $U_q(\overline{\geh})$
module with highest weight $\overline{\eta}$ within 
the integrable highest weight module $V(\lambda)$, 
where $U_q(\overline{\geh})$ stands for the subalgebra of
$\Uq$ generated by $e_i, f_i, t_i ( i \in I\setminus \{ 0 \})$.
\begin{remark}\label{kacwakimoto}
Multiply $q^{-c_j}$ on both sides of (\ref{xbyg}) and 
take $j \rightarrow \infty$ limit.
{}From (\ref{glimit}) and (\ref{xlimit}), the result turns out to be
equivalent with Theorem 3.1 in \cite{KW} when $\mu$ there is dominant
integral.
In this sense (\ref{xbyg}) is a finite $j$ analogue of it.
\end{remark}

One can interpret the 
Kostka-Foulkes polynomial $K_{\xi \mu}(q)$ 
\cite{Ma} 
as a classically restricted 1dsum
for $A^{(1)}_n$.
Consider the level $l$ perfect crystal $B$ corresponding to the 
$l$-fold symmetric tensor representation \cite{KMN2}.
It is parametrized by semistandard tableaux of shape
$(l)$ and entries from $\{0, 1, \ldots, n \}$.
In particular $b(l\Lambda_0)$ is the one with all entries being $n$.
Let $0 \le x_1 \le \cdots \le x_l \le n$ and 
$0 \le y_1 \le \cdots \le y_l \le n$ stand for the 
semistandard tableaux for $b$ and $b' \in B$, respectively.
Then the $H$-function (\ref{defh}) is given by
$H(b \otimes b') = \hbox{ min }_\tau (\sum_{i=1}^l 
\theta(x_i \ge y_{\tau(i)}))$.
Here  
$\theta(\hbox{true}) = 1, \theta(\hbox{false}) = 0$ and the minimum
extends over the degree $l$ symmetric group.
Let $\xi = (\xi_1, \xi_2, \ldots, \xi_{n+1})$ be any partition 
of $lj$ (depth $l(\xi) \le n+1$) and identify it with 
$\sum_{i=1}^n(\xi_i - \xi_{i+1})\overline{\Lambda}_i
\in \overline{P}_{cl}^+$.
Then we have
\begin{equation}
K_{\xi (l^j)}(q) = q^{-lj}
\overline{X}_j(b(l\Lambda_0),0,\xi).
\label{kostka}
\end{equation}
This is just an interpretation of a special case 
of a theorem in \cite{NY} via the classically restricted 
1dsums.
See also \cite{KMOTU2} for another 
extension.
Our picture can be summarized roughly in the following table.

\begin{table}[h]
\caption{}
\begin{center}
\renewcommand{\arraystretch}{2.0}
\begin{tabular}{| c | ccc |}
\hline
{1dsum} & {$g_j$} & {$\overline{X}_j$} & {$X_j$} \\
path & unrestricted & classically restricted & restricted \\
$j < \infty$ & 
${{\hbox{string function of}} \atop {\hbox{Demazure module}}}$ 
 &$\geh$-Kostka & 
restricted $\geh$-Kostka \\
$j \rightarrow \infty$ & string function & 
${\goth{g}_l}/{\overline{\goth{g}}}$ &
$({\goth{g}_{l'} \oplus \goth{g}_l})/{\goth{g}_{l+l'}}$ \\
\hline
\end{tabular}
\end{center}
\end{table}

\noindent
Here $\goth{g}_l$ denotes the affine Lie algebra $\goth{g}$ at level $l$.
By $\geh$-Kostka we generally mean the branching coefficients
of the irreducible $U_q(\overline{\geh})$-modules in the Demazure
module $V_{w^{(jd)}}(\sigma^{-j}(l\Lambda_0))$.
See \cite{KMOTU1}, \cite{KMOTU2} for the 
$U_q(\overline{\geh})$ invariance of the Demazure modules.
\begin{remark}\label{kostkadifference}
Combining (\ref{kostka}) and (\ref{xbarbyg}) one can express 
$K_{\xi (l^j)}(q)$ as an alternating sum over $\overline{W}$.
However the resulting formula is different from the one
on p244 in \cite{Ma}.
\end{remark}

\section{$q$-multinomial formula for $g_j(b,\mu)$}

In this section we present explict formulae for the 
1dsums $g_j(b,\mu)$ in terms of $q$-multinomial coefficients.
We shall also attach the data 
$B, d, i^{(j)}_a, B^{(j)}_a$, etc from \cite{KMOTU1}, which satisfy 
({\bf I}) - ({\bf IV}) in section 1.
Combined with Proposition \ref{pr:by1dsum} or (\ref{simplecase})
they yield a character formula for the Demazure module
$V_{w^{(k)}}(\lambda)$.
We shall only consider level 1 cases of 
$U_{q}(\geh)$ with $\geh$ being classical types: $A^{(1)}_{n}$, 
$B^{(1)}_{n}$, $D^{(1)}_{n}$, $A^{(2)}_{2n-1}$, $A^{(2)}_{2n}$ and 
$D^{(2)}_{n+1}$. 
Other cases, especially higher level cases will be treated elsewhere.
Except the $D^{(2)}_{n+1}$ case, 
the $q$-multinomial formulae for $g_j(b,\mu)$ have been
effectively known in earlier works 
\cite{JMO}, \cite{DJKMO1}, \cite{Kun} on solvable lattice models.
They can be proved by establishing the recursion relation
(\ref{rec1dsum}).

Given a crystal $B$ and an integer vector
with $\sharp B$-components $\gamma = (\gamma_b)_{b \in B}$ 
we shall employ the notations

\begin{eqnarray*}
&\left[ \begin{array}{c} j \\ \gamma \end{array} \right]_q &=
\left \{\begin{array}{cl}
{(q)_j \over \prod_{b \in B}(q)_{\gamma_b}} & \quad 
\mbox{if \, $j = \sum_{b \in B} \gamma_b$ \, 
and $\gamma_b \in {\bf Z}_{\ge 0}$}\\
0& \quad \mbox{otherwise}
\end{array}\right. \\
&(q)_m &= \prod_{i=1}^m(1-q^i) \quad {\rm for }\,  m \in {\bf Z}_{\ge 0}.
\end{eqnarray*}

We shall also use 
$$
\epsilon_s = \left\{
 \begin{array}{ll}
0& {\rm if} \; s \mbox{ is even } \\
1& {\rm if} \; s \mbox{ is odd }. \\
\end{array}\right.
$$
In view of (\ref{1dsumproperty}) we shall assume 
$\mu$ is a level 0 integral weight, i.e., 
$\mu \in P_{cl},\, \langle \mu, c \rangle = 0$
in the rest of the paper.

\subsection{$(A^{(1)}_n, B(\Lambda_1))$ case}

The level 1 perfect crystal 
$B = B(\Lambda_1) = \{0, 1, \ldots, n \}$ can be depicted in the 
crystal graph
 
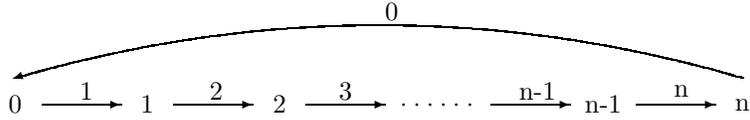
\begin{picture}(300,80)(0,0)
\put(10,10){
\bezier{500}(10,20)(147.5,60)(285,20)
\put(10,20){\vector(-3,-1){1}}
\multiput(20,10)(50,0){3}{\vector(1,0){30}}
\put(5,5){\makebox(10,10){0}}
\put(55,5){\makebox(10,10){1}}
\put(105,5){\makebox(10,10){2}}
\put(32,10){\makebox(10,10){1}}
\put(81,10){\makebox(10,10){2}}
\put(130,10){\makebox(10,10){3}}
\multiput(157,10)(5,0){6}{\circle*{1}}
\multiput(190,10)(55,0){2}{\vector(1,0){30}}
\put(225,5){\makebox(15,10){n-1}}
\put(280,5){\makebox(10,10){n}}
\put(200,10){\makebox(15,10){n-1}}
\put(257,10){\makebox(10,10){n}}
\put(147.5,40){\makebox(10,10){0}}
}
\end{picture}

Elements of $B$ have the weights
$$
wt(b) = \Lambda_{\overline{b+1}} - \La_b\,\, \mbox{ for } b \in B,
$$
where $\overline{x}$ is uniquely specified from $x$ by
$\overline{x} \equiv x$ mod $n+1$ and $0 \le \overline{x} \le n$.
The energy function is given by
$$
H(b\otimes b') = \left\{\begin{array}{cc}
0& {\rm if} \; b <  b' \\
1& {\rm if} \; b \ge b'.
\end{array}\right.
$$
Due to the Dynkin diagram symmetry it suffices to consider
the case $\la = \Lambda_0$.
Then we have the result \cite{KMOTU1}:
\begin{eqnarray*}
&&d = n,\quad \lambda_j = \Lambda_{\overline{ -j}},\quad
\bbar_j = \overline{ -j},  \\
&&B^{(j)}_a = \{ \overline{-j}, \overline{-j+1}, \ldots, \overline{-j+a} \}
\quad 1 \le a \le n, \\
&&i^{(j)}_a =  \overline{-j+a}. \\
\end{eqnarray*}
\begin{proposition}[{\rm cf. \cite{JMO}, \cite{DJKMO1}}]
For $j \in {\bf Z}_{\ge 0}, \, b \in B$ and 
$\mu = (\mu_i)_{i \in B} = 
(\mu_n - \mu_0)\Lambda_0 + (\mu_0 - \mu_1) \Lambda_1 + \cdots +
(\mu_{n-1} - \mu_n)\Lambda_n \in P_{cl}$, we have
$$
g_j(b,\mu) = q^{\frac{1}{2} \sum_{i\in B}\mu_{i}(\mu_{i}-1)+
\sum_{i\in B}H(b \ot i)\mu_i}{j \brack \mu}_q.
$$
\end{proposition}

\subsection{$(B^{(1)}_n, B(\Lambda_1))$ Case}

The level 1 perfect crystal 
$B = B(\Lambda_1) =
\{1,2, \ldots, n, 0, \overline{n}, \ldots, \overline{1} \}$ 
is depicted by the crystal graph

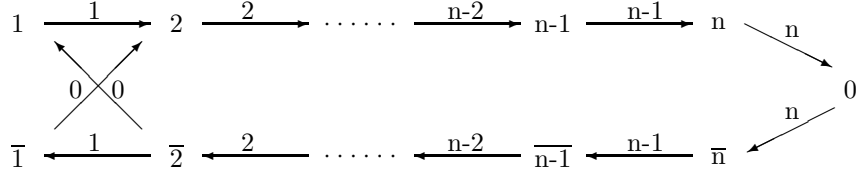
\begin{picture}(300,100)(0,0)
\put(0,15){
\multiput(20,60)(60,0){2}{\vector(1,0){40}}
\put(5,55){\makebox(10,10){1}}
\put(65,55){\makebox(10,10){2}}
\put(34,60){\makebox(10,10){1}}
\put(92,60){\makebox(10,10){2}}
\multiput(127,60)(5,0){6}{\circle*{1}}
\multiput(160,60)(65,0){2}{\vector(1,0){40}}
\put(205,55){\makebox(15,10){n-1}}
\put(270,55){\makebox(10,10){n}}
\put(172,60){\makebox(15,10){n-2}}
\put(240,60){\makebox(15,10){n-1}}
\multiput(60,10)(60,0){2}{\vector(-1,0){40}}
\put(5,5){\makebox(10,10){$\overline{1}$}}
\put(65,5){\makebox(10,10){$\overline{2}$}}
\put(34,10){\makebox(10,10){1}}
\put(92,10){\makebox(10,10){2}}
\multiput(127,10)(5,0){6}{\circle*{1}}
\multiput(200,10)(65,0){2}{\vector(-1,0){40}}
\put(205,5){\makebox(15,10){$\overline{\mbox{n-1}}$}}
\put(270,5){\makebox(10,10){$\overline{\mbox{n}}$}}
\put(172,10){\makebox(15,10){n-2}}
\put(240,10){\makebox(15,10){n-1}}
\put(320,30){\makebox(10,10){0}}
\put(285,60){\vector(2,-1){33}}
\put(319,28){\vector(-2,-1){33}}
\put(298,52){\makebox(10,10){n}}
\put(298,21){\makebox(10,10){n}}
\put(57,20){\vector(-1,1){33}}
\put(24,20){\vector(1,1){33}}
\put(27,30){\makebox(10,10){0}}
\put(43,30){\makebox(10,10){0}}
}
\end{picture}

Elements of $B$ have the weights
\begin{eqnarray*}
&&wt(b) = - wt(\overline{b}) = \left\{
\begin{array}{ll}
\Lambda_b - \Lambda_{b-1} & b = 1 \mbox{ or } 3 \le b \le n-1 \\
\La_2 - \La_1 - \La_0 & b = 2 \\
2\La_n - \La_{n-1} & b = n \\
\end{array}\right. \\
&&wt(0) = 0.\\
\end{eqnarray*}
We introduce an order $\prec$ on $B$ by
$$
1 \prec \cdots  \prec n \prec 0 \prec \overline{n} \prec \cdots \prec 
\overline{1}.
$$
This and similar $\prec$ will be used in the subsequent subsections
just for convenience and should not be confused with the Bruhat order.
The energy function is given by 
$$
H(b\otimes b') = \left\{\begin{array}{cc}
0& {\rm if} \; b \prec b' \\
1& {\rm if} \; b \succeq b'.
\end{array}\right.
$$
with the exceptions:
$$
H(0\otimes 0) =0, \quad  H(1\otimes \overline{1}) = -1.
$$
There are 3 level 1 dominant integral weights
$(P^+_{cl})_1 = \{ \Lambda_0, \Lambda_1, \Lambda_n \}$.
Due to the Dynkin diagram symmetry it suffices to consider
$\la = \Lambda_0$ and $\Lambda_n$.
In both cases we have $d = 2n-1$.
The other data given in \cite{KMOTU1} reads 
\begin{eqnarray*}
&\la = \Lambda_0 &{\rm case}:\\
&&\la_j = \La_{\epsilon_j},
\quad 
\bbar_j = \left\{
        \begin{array}{ll}
1& \; j: {\rm even} \\
\overline{1}& \; j: {\rm odd} \\
\end{array}\right.,\\
&&B^{(j)}_a = \left\{
\begin{array}{ll}
\{\bbar_j, 2, \ldots, a+1 \} & 1 \le a \le n-1\\
\{\bbar_j, 2, \ldots, n, 0, \overline{n}, \overline{n-1}, \ldots, 
\overline{2n-a} \} & n \le a \le 2n-2\\
\end{array}\right., \\
&&i^{(j)}_a = \left\{
\begin{array}{ll}
\epsilon_{1-j} & a=1 \mbox{ or } a = 2n-1\\
{\rm min}(a,2n-a)& 2 \le a \le 2n-2.\\
\end{array}\right.
\end{eqnarray*}
\begin{eqnarray*}
&\la = \Lambda_n &{\rm case}:\\
&&\la_j = \La_n,
\quad 
\bbar_j = 0,\\
&&B^{(j)}_a = \left\{
\begin{array}{ll}
\{0, \overline{n}, \overline{n-1}, \ldots, \overline{n+1-a} \} & 1 \le a \le n-1\\
\{0, \overline{n}, \ldots, \overline{2}, 1 \} & a = n\\
\{0, \overline{n}, \overline{n-1}, \ldots, \overline{1}, 1, 2, \ldots, a-n+1 \} 
& n+1 \le a \le 2n-1\\
\end{array}\right., \\
&&i^{(j)}_a = \left\{
        \begin{array}{ll}
n+1-a & 1 \le a \le n-1 \\
a-n & n \le a \le 2n-1. \\
\end{array}\right.
\end{eqnarray*}
For $\la = \La_n$ one can also make another 
choice of $B^{(j)}_a$ and $i^{(j)}_a$ as
\begin{eqnarray*}
&&B^{(j)}_a = \left\{
\begin{array}{ll}
\{0, \overline{n}, \overline{n-1}, \ldots, \overline{n+1-a} \} & 0 \le a \le n\\
\{0, \overline{n}, \overline{n-1}, \ldots, \overline{1}, 1, 2, \ldots, a-n+1 \} 
& n+1 \le a \le 2n-1\\
\end{array}\right., \\
&&i^{(j)}_a = \left\{
        \begin{array}{ll}
n+1-a & 1 \le a \le n+1 \\
a-n & n+2 \le a \le 2n-1. \\
\end{array}\right.
\end{eqnarray*}
In any case, $B^{(j)}_0 = \{ \bbar_j \}$ and $B^{(j)}_d = B$ hold.
Let us parametrize the level 0 elements $\mu \in P_{cl}$ by
$(\mu_i)_{i=1}^n \in {\bf Z}^n$ as
\[
\mu = (- \mu_1 - \mu_2)\La_0 + (\mu_1 - \mu_2)\La_1 + \cdots +
(\mu_{n-1} - \mu_n) \Lambda_{n-1} + 2\mu_n \La_n.
\]
\begin{proposition}[{\rm cf. \cite{DJKMO1}}]
For $j \in {\bf Z}_{\ge 0}, b \in B$ and the above $\mu \in P_{cl}$, we have
\begin{eqnarray*}
g_j(b,\mu) &=& \mathop{{\sum}^{*}}_\gamma q^{\frac{1}{2} 
\sum_{i\in B}\gamma_{i}(\gamma_{i}-1)-\gamma_s\gamma_{\overline{s}}+
\sum_{i \in B}H(b \ot i)\gamma_{i}}{j\brack \gamma}_q,\\
\end{eqnarray*}
where $s = n$ if $b \succ 0$ and $s = 1$ if $b \prec 0$.
When $b = 0$, either choice $s=n$ or $s=1$ is valid.
The sum ${\sum}^{*}_\gamma$ extends over 
$\gamma = (\gamma_i)_{i \in B}  \in ({\bf Z}_{\ge 0})^{2n+1}$
such that
$$
\gamma_{i} - \gamma_{\overline{i}} = \mu_{i}\;\; 
for\; i = 1, \cdots, n, \quad \sum_{i\in B} \gamma_{i} = j.
$$
\end{proposition}

\subsection{$(D^{(1)}_{n}, B(\Lambda_{1}))$ case}

The level 1 perfect crystal 
$B = B(\Lambda_1) =
\{1,2, \ldots, n, \overline{n}, \ldots, \overline{1} \}$ 
is depicted by the crystal graph
%

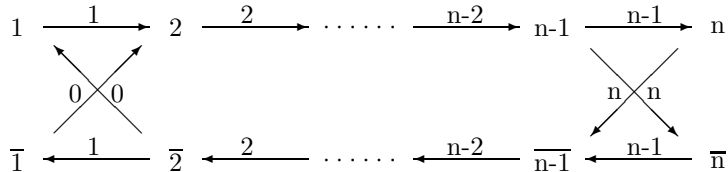
\begin{picture}(300,100)(0,0)
\put(10,15){
\multiput(20,60)(60,0){2}{\vector(1,0){40}}
\put(5,55){\makebox(10,10){1}}
\put(65,55){\makebox(10,10){2}}
\put(34,60){\makebox(10,10){1}}
\put(92,60){\makebox(10,10){2}}
\multiput(127,60)(5,0){6}{\circle*{1}}
\multiput(160,60)(65,0){2}{\vector(1,0){40}}
\put(205,55){\makebox(15,10){n-1}}
\put(270,55){\makebox(10,10){n}}
\put(172,60){\makebox(15,10){n-2}}
\put(240,60){\makebox(15,10){n-1}}
\multiput(60,10)(60,0){2}{\vector(-1,0){40}}
\put(5,5){\makebox(10,10){$\overline{1}$}}
\put(65,5){\makebox(10,10){$\overline{2}$}}
\put(34,10){\makebox(10,10){1}}
\put(92,10){\makebox(10,10){2}}
\multiput(127,10)(5,0){6}{\circle*{1}}
\multiput(200,10)(65,0){2}{\vector(-1,0){40}}
\put(205,5){\makebox(15,10){$\overline{\mbox{n-1}}$}}
\put(270,5){\makebox(10,10){$\overline{\mbox{n}}$}}
\put(172,10){\makebox(15,10){n-2}}
\put(240,10){\makebox(15,10){n-1}}
\put(260,52){\vector(-1,-1){33}}
\put(227,52){\vector(1,-1){33}}
\put(246,30){\makebox(10,10){n}}
\put(231,30){\makebox(10,10){n}}
\put(57,20){\vector(-1,1){33}}
\put(24,20){\vector(1,1){33}}
\put(27,30){\makebox(10,10){0}}
\put(43,30){\makebox(10,10){0}}
}
\end{picture}

Elements of $B$ have the weights
$$
wt(b) = - wt(\overline{b}) = \left\{
\begin{array}{ll}
\Lambda_b - \Lambda_{b-1} & b \ne 2, n-1 \\
\La_2 - \La_1 - \La_0 & b = 2 \\
\La_n + \La_{n-1} - \La_{n-2} & b = n-1. \\
\end{array}\right. 
$$
We introduce an order on $B$ by
$$
1 \prec \cdots  \prec n-1 \prec {n \atop\overline{n}} \prec 
\overline{n-1} \prec \cdots \prec \overline{1}.
$$
There is no order between $n$ and $\overline{n}$.
The energy function is given by
$$
H(b\otimes b') = \left\{\begin{array}{cc}
0& {\rm if} \; b \prec b', \\
1& {\rm if} \; b \succeq b'.
\end{array}\right.
$$
with the exceptions:
$$
H(n \otimes \overline{n}) = H(\overline{n}\otimes n) =0, \quad  
H(1\otimes \overline{1}) = -1.
$$
There are 4 level 1 dominant integral weights
$(P^+_{cl})_1 = \{ \Lambda_0, \Lambda_1, \La_{n-1}, \Lambda_n \}$.
Due to the Dynkin diagram symmetry it suffices to consider
$\la = \Lambda_0$.
Then the result in \cite{KMOTU1} reads
\begin{eqnarray*}
&&d  = 2n-2,\quad \la_j = \La_{\epsilon_j},
\quad 
\bbar_j = \left\{
        \begin{array}{ll}
1& \; j: {\rm even} \\
\overline{1}& \; j: {\rm odd} \\
\end{array}\right.,\\
&&B^{(j)}_a = \left\{
\begin{array}{ll}
\{\bbar_j, 2, \ldots, a+1 \} & 1 \le a \le n-2\\
\{\bbar_j, 2, \ldots, n-1, \overline{n} \} & a = n-1\\
\{\bbar_j, 2, \ldots, n, \overline{n}, \overline{n-1}, \ldots, 
\overline{2n-1-a} \} & n \le a \le 2n-3\\
\end{array}\right., \\
&&i^{(j)}_a = \left\{
\begin{array}{ll}
\epsilon_{1-j} & a=1, 2n-2\\
{\rm min}(a,2n-1-a)& a \neq 1, n-1, 2n-2\\
2n-1-a & a=n-1.\\
\end{array}\right.
\end{eqnarray*}
One can also make another choice of 
$B^{(j)}_a$ and $i^{(j)}_a$ as
\begin{eqnarray*}
&&B^{(j)}_a = \left\{
\begin{array}{ll}
\{\bbar_j, 2, \ldots, a+1 \} & 1 \le a \le n-1\\
\{\bbar_j, 2, \ldots, n, \overline{n}, \overline{n-1}, \ldots, 
\overline{2n-1-a} \} & n \le a \le 2n-3\\
\end{array}\right., \\
&&i^{(j)}_a = \left\{
\begin{array}{ll}
\epsilon_{1-j} & a=1, 2n-2\\
{\rm min}(a,2n-1-a)& a \neq 1, n, 2n-2\\
n & a=n.\\
\end{array}\right. 
\end{eqnarray*}
In any case, $B^{(j)}_0 = \{ \bbar_j \}$ and $B^{(j)}_d = B$ hold.
Let us parametrize the level 0 elements $\mu \in P_{cl}$ by
$(\mu_i)_{i=1}^n \in {\bf Z}^n$ as
\[
\mu = (- \mu_1 - \mu_2)\La_0 + (\mu_1 - \mu_2)\La_1 + \cdots +
(\mu_{n-1} - \mu_n) \Lambda_{n-1} + (\mu_{n-1} + \mu_n) \La_n.
\]
\begin{proposition}[{\rm cf. \cite{DJKMO1}}]
For $j \in {\bf Z}_{\ge 0}, b \in B$ and the above $\mu \in P_{cl}$, we have
\begin{eqnarray*}
g_j(b,\mu) &=& \mathop{{\sum}^{*}}_\gamma q^{\frac{1}{2} 
\sum_{i\in B}\gamma_{i}(\gamma_{i}-1)-\gamma_s\gamma_{\overline{s}}+
\sum_{i \in B}H(b \ot i)\gamma_{i}}{j\brack \gamma}_q,\\
\end{eqnarray*}
where $s = n$ if $b \in \{\overline{n}, \ldots, \overline{1} \}$ 
and $s = 1$ if $b \in \{ 1, \ldots, n \}$.
The sum ${\sum}^{*}_\gamma$ extends over 
$\gamma = (\gamma_i)_{i \in B}  \in ({\bf Z}_{\ge 0})^{2n}$
such that
$$
\gamma_{i} - \gamma_{\overline{i}} = \mu_{i}\;\; 
for\; i = 1, \cdots, n, \quad \sum_{i\in B} \gamma_{i} = j.
$$
\end{proposition}
This is a very similar form to the $B^{(1)}_n$ case.

\subsection{$(A^{(2)}_{2n-1}, B(\Lambda_{1}))$ case}

The level 1 perfect crystal 
$B = B(\Lambda_1) =
\{1,2, \ldots, n, \overline{n}, \ldots, \overline{1} \}$ 
is depicted by the crystal graph
%

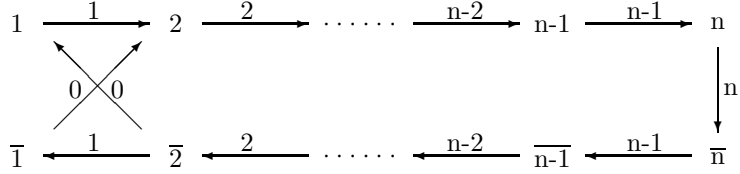
\begin{picture}(300,100)(0,0)
\put(10,15){
\multiput(20,60)(60,0){2}{\vector(1,0){40}}
\put(5,55){\makebox(10,10){1}}
\put(65,55){\makebox(10,10){2}}
\put(34,60){\makebox(10,10){1}}
\put(92,60){\makebox(10,10){2}}
\multiput(127,60)(5,0){6}{\circle*{1}}
\multiput(160,60)(65,0){2}{\vector(1,0){40}}
\put(205,55){\makebox(15,10){n-1}}
\put(270,55){\makebox(10,10){n}}
\put(172,60){\makebox(15,10){n-2}}
\put(240,60){\makebox(15,10){n-1}}
\multiput(60,10)(60,0){2}{\vector(-1,0){40}}
\put(5,5){\makebox(10,10){$\overline{1}$}}
\put(65,5){\makebox(10,10){$\overline{2}$}}
\put(34,10){\makebox(10,10){1}}
\put(92,10){\makebox(10,10){2}}
\multiput(127,10)(5,0){6}{\circle*{1}}
\multiput(200,10)(65,0){2}{\vector(-1,0){40}}
\put(205,5){\makebox(15,10){$\overline{\mbox{n-1}}$}}
\put(270,5){\makebox(10,10){$\overline{\mbox{n}}$}}
\put(172,10){\makebox(15,10){n-2}}
\put(240,10){\makebox(15,10){n-1}}
\put(275,51){\vector(0,-1){32}}
\put(275,30){\makebox(10,10){n}}
\put(57,20){\vector(-1,1){33}}
\put(24,20){\vector(1,1){33}}
\put(27,30){\makebox(10,10){0}}
\put(43,30){\makebox(10,10){0}}
}
\end{picture}

Elements of $B$ have the weights
$$
wt(b) = - wt(\overline{b}) = \left\{
\begin{array}{ll}
\Lambda_b - \Lambda_{b-1} & b \ne 2\\
\La_2 - \La_1 - \La_0 & b = 2. \\
\end{array}\right. 
$$
We introduce an order on $B$ by
$$
1 \prec \cdots  \prec n \prec \overline{n} \prec \cdots \prec 
\overline{1}.
$$
The energy function is given by 
$$
H(b\otimes b') = \left\{\begin{array}{cc}
0& {\rm if} \; b \prec b' \\
1& {\rm if} \; b \succeq b',
\end{array}\right.
$$
with the exception:
$$
H(1\otimes \overline{1}) = -1.
$$
There are 2 level 1 dominant integral weights
$(P^+_{cl})_1 = \{ \Lambda_0, \Lambda_1 \}$.
Due to the Dynkin diagram symmetry it suffices to consider
$\la = \Lambda_0$.
Then the result in \cite{KMOTU1} reads
\begin{eqnarray*}
&&d  = 2n-1,\quad \la_j = \La_{\epsilon_j},
\quad 
\bbar_j = \left\{
        \begin{array}{ll}
1& \; j: {\rm even} \\
\overline{1}& \; j: {\rm odd} \\
\end{array}\right.,\\
&&B^{(j)}_a = \left\{
\begin{array}{ll}
\{\bbar_j, 2, \ldots, a+1 \} & 1 \le a \le n-1\\
\{\bbar_j, 2, \ldots, n, \overline{n}, \overline{n-1}, \ldots, 
\overline{2n-a} \} & n \le a \le 2n-2\\
\end{array}\right., \\
&&i^{(j)}_a = \left\{
\begin{array}{ll}
\epsilon_{1-j} & a=1, 2n-1\\
{\rm min}(a,2n-a)& 2 \le a \le 2n-2.\\
\end{array}\right.
\end{eqnarray*}
$B^{(j)}_0 = \{ \bbar_j \}$ and $B^{(j)}_d = B$ hold.
Let us parametrize the level 0 elements $\mu \in P_{cl}$ by
$(\mu_i)_{i=1}^n \in {\bf Z}^n$ as
\[
\mu = (- \mu_1 - \mu_2)\La_0 + (\mu_1 - \mu_2)\La_1 + \cdots +
(\mu_{n-1} - \mu_n) \Lambda_{n-1} + \mu_n \La_n.
\]
\begin{proposition}[{\rm cf. \cite{Kun}}]
For $j \in {\bf Z}_{\ge 0}, b \in B$ and the above $\mu \in P_{cl}$, we have
\begin{eqnarray*}
g_j(b,\mu) &=& \mathop{{\sum}^{*}}_\gamma q^{Q}
\frac{(q^2)_{\gamma_1 + \gamma_{\overline{1}}}(q)_j}
{(q^2)_{\gamma_1}(q^2)_{\gamma_{\overline{1}}}
(q)_{\gamma_1 + \gamma_{\overline{1}}}
\prod_{i=2}^n(q)_{\gamma_i}(q)_{\gamma_{\overline{i}}}}
G(b,\mu,\gamma),\\
Q &=& \frac{1}{2} \sum_{i\in B}\gamma_{i}(\gamma_{i}-1)-
\gamma_{1}\gamma_{\overline{1}}+
\sum_{i\in B}H(b \ot i)\gamma_{i}, \\
G(b,\mu,\gamma) &=& \left\{ \begin{array}{cc}
1 & if\;\; b = 1\;or\;\overline{1} \\
\frac{ q^{\frac{1}{2}\mu_{1}} + q^{-\frac{1}{2}\mu_{1}} }
 { q^{\frac{1}{2}(\gamma_{1}+\gamma_{\overline{1}})} + q^{-\frac{1}{2}
 (\gamma_{1}+\gamma_{\overline{1}}) }} & otherwise.
\end{array}\right.,
\end{eqnarray*}
The sum ${\sum}^{*}_\gamma$ extends over 
$\gamma = (\gamma_i)_{i \in B}  \in ({\bf Z}_{\ge 0})^{2n}$
such that
$$
\gamma_{i} - \gamma_{\overline{i}} = \mu_{i}\;\; 
for\; i = 1, \cdots, n, \quad \sum_{i\in B} \gamma_{i} = j.
$$
\end{proposition}

\subsection{$(A^{(2)}_{2n}, B(0) \oplus B(\Lambda_{1}))$ case}

For a technical reason, we take the opposite ordering for the labeling of 
vertices of the Dynkin diagram from \cite{KMOTU1}.
The level 1 perfect crystal 
$B = B(0) \oplus B(\Lambda_1) =
\{1,2, \ldots, n, 0, \overline{n}, \ldots, \overline{1} \}$ 
 is depicted by the crystal graph

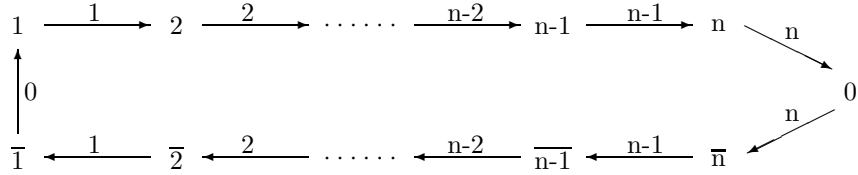
\begin{picture}(300,100)(0,0)
\put(0,15){
\multiput(20,60)(60,0){2}{\vector(1,0){40}}
\put(5,55){\makebox(10,10){1}}
\put(65,55){\makebox(10,10){2}}
\put(34,60){\makebox(10,10){1}}
\put(92,60){\makebox(10,10){2}}
\multiput(127,60)(5,0){6}{\circle*{1}}
\multiput(160,60)(65,0){2}{\vector(1,0){40}}
\put(205,55){\makebox(15,10){n-1}}
\put(270,55){\makebox(10,10){n}}
\put(172,60){\makebox(15,10){n-2}}
\put(240,60){\makebox(15,10){n-1}}
\multiput(60,10)(60,0){2}{\vector(-1,0){40}}
\put(5,5){\makebox(10,10){$\overline{1}$}}
\put(65,5){\makebox(10,10){$\overline{2}$}}
\put(34,10){\makebox(10,10){1}}
\put(92,10){\makebox(10,10){2}}
\multiput(127,10)(5,0){6}{\circle*{1}}
\multiput(200,10)(65,0){2}{\vector(-1,0){40}}
\put(205,5){\makebox(15,10){$\overline{\mbox{n-1}}$}}
\put(270,5){\makebox(10,10){$\overline{\mbox{n}}$}}
\put(172,10){\makebox(15,10){n-2}}
\put(240,10){\makebox(15,10){n-1}}
\put(320,30){\makebox(10,10){0}}
\put(285,60){\vector(2,-1){33}}
\put(319,28){\vector(-2,-1){33}}
\put(298,52){\makebox(10,10){n}}
\put(298,21){\makebox(10,10){n}}
\put(10,19){\vector(0,1){32}}
\put(10,30){\makebox(10,10){0}}
}
\end{picture}

Elements of $B$ have the weights
\begin{eqnarray*}
&&wt(b) = - wt(\overline{b}) = \left\{
\begin{array}{ll}
\Lambda_b - \Lambda_{b-1} & 1 \le b \le n-1 \\
2\La_n - \La_{n-1} & b = n \\
\end{array}\right. \\
&&wt(0) = 0.\\
\end{eqnarray*}
We introduce an order $\prec$ on $B$ by
$$
1 \prec \cdots  \prec n \prec 0 \prec \overline{n} \prec \cdots \prec 
\overline{1}.
$$
The energy function is given by 
$$
H(b\otimes b') = \left\{\begin{array}{cc}
0& if \; b \prec b' \\
1& if \; b \succeq b',
\end{array}\right.
$$
with the exception:
$$
H(0\otimes 0) =0.
$$
There is a unique level 1 dominant integral weight
$(P^+_{cl})_1 = \{ \Lambda_n \}$.
Thus we set $\la = \La_n$, for which 
the result in \cite{KMOTU1} reads
\begin{eqnarray*}
&&d  = 2n,\quad \la_j = \La_n,
\quad 
\bbar_j = 0,\\
&&B^{(j)}_a = \left\{
\begin{array}{ll}
\{0, \overline{n}, \overline{n-1}, \ldots, 
\overline{n+1-a} \} & 1 \le a \le n\\
\{0, \overline{n}, \overline{n-1}, \ldots, \overline{1},
1, 2, \ldots, a-n \} & n+1 \le a \le 2n\\
\end{array}\right., \\
&&i^{(j)}_a = \vert n+1-a \vert \quad 1 \le a \le 2n.
\end{eqnarray*}
Let us parametrize the level 0 elements $\mu \in P_{cl}$ by
$(\mu_i)_{i=1}^n \in {\bf Z}^n$ as
\[
\mu = - \mu_1 \La_0 + (\mu_1 - \mu_2)\La_1 + \cdots +
(\mu_{n-1} - \mu_n) \Lambda_{n-1} + 2\mu_n \La_n.
\]
\begin{proposition}[{\rm cf. \cite{Kun}}]
For $j \in {\bf Z}_{\ge 0}, b \in B$ and the above $\mu \in P_{cl}$, we have
\begin{eqnarray*}
g_j(b,\mu) &=& \mathop{{\sum}^{*}}_\gamma q^{\frac{1}{2} 
\sum_{i\in B, i \neq 0}\gamma_{i}(\gamma_{i}-1)+
\sum_{i \in B}H(b \ot i)\gamma_{i}}{j\brack \gamma}_q,\\
\end{eqnarray*}
where the sum ${\sum}^{*}_\gamma$ extends over 
$\gamma = (\gamma_i)_{i \in B}  \in ({\bf Z}_{\ge 0})^{2n+1}$
such that
$$
\gamma_{i} - \gamma_{\overline{i}} = \mu_{i}\;\; 
for\; i = 1, \cdots, n, \quad \sum_{i\in B} \gamma_{i} = j.
$$
\end{proposition}

\subsection{$(D^{(2)}_{n+1}, B(0) \oplus B(\Lambda_1))$ Case}

The level 1 perfect crystal 
$B = B(0) \oplus B(\Lambda_1) =
\{1,2, \ldots, n, 0, \overline{n}, \ldots, \overline{1}, \phi \}$ 
is depicted by the crystal graph

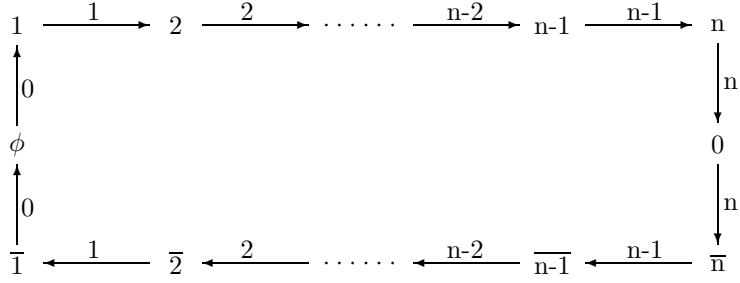
\begin{picture}(300,130)(0,0)
\put(10,10){
\multiput(20,100)(60,0){2}{\vector(1,0){40}}
\put(5,95){\makebox(10,10){1}}
\put(65,95){\makebox(10,10){2}}
\put(34,100){\makebox(10,10){1}}
\put(92,100){\makebox(10,10){2}}
\multiput(127,100)(5,0){6}{\circle*{1}}
\multiput(160,100)(65,0){2}{\vector(1,0){40}}
\put(205,95){\makebox(15,10){n-1}}
\put(270,95){\makebox(10,10){n}}
\put(172,100){\makebox(15,10){n-2}}
\put(240,100){\makebox(15,10){n-1}}
\multiput(60,10)(60,0){2}{\vector(-1,0){40}}
\put(5,5){\makebox(10,10){$\overline{1}$}}
\put(65,5){\makebox(10,10){$\overline{2}$}}
\put(34,10){\makebox(10,10){1}}
\put(92,10){\makebox(10,10){2}}
\multiput(127,10)(5,0){6}{\circle*{1}}
\multiput(200,10)(65,0){2}{\vector(-1,0){40}}
\put(205,5){\makebox(15,10){$\overline{\mbox{n-1}}$}}
\put(270,5){\makebox(10,10){$\overline{\mbox{n}}$}}
\put(172,10){\makebox(15,10){n-2}}
\put(240,10){\makebox(15,10){n-1}}
\put(270,50){\makebox(10,10){0}}
\multiput(275,47)(0,46){2}{\vector(0,-1){30}}
\multiput(275,27)(0,46){2}{\makebox(10,10){n}}
\put(5,50){\makebox(10,10){$\phi$}}
\multiput(10,17)(0,45){2}{\vector(0,1){30}}
\multiput(9,26)(0,45){2}{\makebox(10,10){0}}
}
\end{picture}

Elements of $B$ have the weights
\begin{eqnarray*}
&&wt(b) = - wt(\overline{b}) = \left\{
\begin{array}{ll}
\La_1-2\La_0 & b = 1 \\
\Lambda_b - \Lambda_{b-1} & 2 \le b \le n-1 \\
2\La_n - \La_{n-1} & b = n \\
\end{array}\right. \\
&&wt(0) = wt(\phi) = 0.\\
\end{eqnarray*}
We introduce an order $\prec$ on $B\setminus \{ \phi \}$ by
$$
1 \prec \cdots  \prec n \prec 0 \prec \overline{n} \prec \cdots \prec 
\overline{1}.
$$
The energy function is given by 
$$
H(b\otimes b') = \left\{\begin{array}{ll}
0& {\rm if} \; b \prec b' \mbox{ or } (b, b') = (0,0), (\phi,\phi) \\
1& \mbox{if one and only one of } b \mbox{ and } b'\mbox{ is } \phi \\
2& {\rm if} \; b \succeq b' \mbox{ and }
(b, b') \neq (0,0), (\phi,\phi). \\
\end{array}\right.
$$
There are 2 level 1 dominant integral weights
$(P^+_{cl})_1 = \{ \Lambda_0, \Lambda_n \}$.
Due to the Dynkin diagram symmetry it suffices to consider
$\la = \Lambda_0$.
Then the result in \cite{KMOTU1} reads
\begin{eqnarray*}
&&d  = 2n,\quad \la_j = \La_0,
\quad 
\bbar_j = \phi,\\
&&B^{(j)}_a = \left\{
\begin{array}{ll}
\{\phi, 1, 2, \ldots, a \} & 0 \le a \le n \\
\{\phi, 1, \ldots, n, 0, \overline{n}, \overline{n-1}, \ldots, 
\overline{2n+1-a} \} & n+1 \le a \le 2n \\
\end{array}\right. , \\
&&i^{(j)}_a = {\rm min}(a-1,2n+1-a) \quad 1 \le a \le 2n.\\
\end{eqnarray*}
Let us parametrize the level 0 elements $\mu \in P_{cl}$ by
$(\mu_i)_{i=1}^n \in {\bf Z}^n$ as
\[
\mu = - 2\mu_1 \La_0 + (\mu_1 - \mu_2)\La_1 + \cdots +
(\mu_{n-1} - \mu_n) \Lambda_{n-1} + 2\mu_n \La_n.
\]
\begin{proposition}
For $j \in {\bf Z}_{\ge 0}, b \in B$ and the above $\mu \in P_{cl}$, we have
\begin{eqnarray*}
g_j(b,\mu) &=& \mathop{{\sum}^{*}}_\gamma 
q^{\sum_{i\in B, i \neq 0, \phi}\gamma_{i}(\gamma_{i}-1)+
\sum_{i \in B}H(b \ot i)\gamma_{i}}{j\brack \gamma}_{q^2},\\
\end{eqnarray*}
where the sum ${\sum}^{*}_\gamma$ extends over 
$\gamma = (\gamma_i)_{i \in B}  \in ({\bf Z}_{\ge 0})^{2n+2}$
such that
$$
\gamma_{i} - \gamma_{\overline{i}} = \mu_{i}\;\; 
for\; i = 1, \cdots, n, \quad \sum_{i\in B} \gamma_{i} = j.
$$
\end{proposition}

\section{Discussion}

We have shown that 
various characters can be viewed in a unified way as the 1dsums 
under the path realization of Demazure crystals.
Our picture is summarized in Table 1 in the end of section 2.
It is yet another task to actually evaluate these 1dsums.
In this paper it has been done in section 3 for 
level 1 cases of the unrestricted 1dsum $g_j$.
Substitution of them into Proposition \ref{pr:xbyg}
generates formulae also for $X_j$ and $\overline{X}_j$.
As seen explicitly there, the results necessarily involve 
alternating signs from the Weyl group signature.
Such formulae are sometimes called bosonic.
In this respect it is interesting also 
to seek fermionic formulae.
By this one roughly means those series or polynomials
which are free of signs, admit a quasi-particle interpretation
or have an origin in string hypotheses in the Bethe ansatz, etc.
Formulae with such features have been 
explored extensively for several cases of 
$\overline{X}_j$ and $X_j$ in our Table 1 by many authors.
See for example \cite{Ki} and references therein.
On the other hand relatively fewer fermionic formulae seem
known or even conjectured for $g_j$.
A possible reason for this is that $g_j$ does {\em not}
correspond to a counting of highest weight vectors 
as opposed to $\overline{X}_j$ and $X_j$.
We hope to discuss this point further and higher level cases
in near future.

{\it Acknowledgement}.
A.K. would like to thank 
I. Cherednik, T. Nakanishi and A. Varchenko
for their warm hospitality at University of North Carolina.
K.C.M. is supported in part by NSA/MSP Grant No. 96-1-0013.
K.C.M. and M.O. thank O. Foda for his hospitality at
the University of Melbourne, and for collaboration in their
earlier work \cite{FMO}.
M.O. would like to thank M. Wakimoto for discussions
and sending the paper \cite{KW}.

\end{document}